\documentclass[11pt, a4paper]{article}
\usepackage{t1enc}
\usepackage[latin1]{inputenc}
\usepackage[english]{babel}
\usepackage{amsmath,amsthm}
\usepackage{amsfonts}
\usepackage{latexsym}
\usepackage[dvips]{graphicx}
\usepackage{float}
\usepackage[natural]{xcolor}
\usepackage{algorithm}
\usepackage{algorithmic}
\usepackage{enumerate,enumitem}
\usepackage{multirow}
\usepackage{xcolor}
\usepackage[colorlinks,linkcolor=blue]{hyperref}
\usepackage{lineno}
\usepackage{setspace}
\usepackage{fullpage}
\usepackage[title]{appendix}
\usepackage{todonotes}
\usepackage{amssymb}

\usepackage{listings}
\usepackage{bbm}

\usepackage{fullpage}
\usepackage{hyperref}
\usepackage[margin=2.3cm]{geometry}
\usepackage{enumitem,verbatim}

\newtheorem{theorem}{Theorem}[section]
\newtheorem{claim}[theorem]{Claim}
\newtheorem{example}[theorem]{Example}

\newtheorem{lemma}[theorem]{Lemma}

\theoremstyle{definition}

\newtheorem{definition}[theorem]{Definition}

\title{\LARGE Large cliques in graphs with forbidden semi-induced structures}

\author{Nannan Chen \thanks{School of Mathematics, Shandong University, Jinan, China, and Extremal Combinatorics and Probability Group (ECOPRO), Institute for Basic Science (IBS), Daejeon, South Korea. Email: {\tt chennannan@sdu.edu.cn}. Suppported by the China Scholarship Council and IBS-R029-C4.}
\and
Yulai Ma \thanks{Center for Combinatorics and LPMC, Nankai University, Tianjin 300071, P.R.China. Email: {\tt ylma92@163.com}.}
\and
Fan Yang \thanks{Data Science Institute, Shandong University, Jinan, Shandong, China, and Extremal Combinatorics and Prob- ability Group (ECOPRO), Institute for Basic Science (IBS), Daejeon, South Korea. Email: {\tt fyang@sdu.edu.cn}. Supported by Natural Science Foundation of China (12301447,12571367), by Natural Science Foundation of Shandong Province (ZR2024QA056), by the China Postdoctoral Science Foundation (12570073310023) and by China Scholarship Council and IBS-R029-C4.}
}
\date{}

\begin{document}
\maketitle

\begin{abstract}
In 2022, Holmsen showed that any graph with at least \( c \binom{n}{r} \) \(r\)-cliques but no induced complete $r$-partite graph $K_{2,\ldots, 2}$ must contain a clique of order \(\Omega(c^{2^{r-1}} n)\). 
In this paper, we study graphs forbidding semi-induced substructures and show that every $n$-vertex graph $G$ containing at least $c\binom{n}{r}$  copies of $K_r$ (for some constant $c>0$) and forbidding semi-induced substructures, related to $K_{2,\ldots, 2}$, must contain a clique of order $\Omega(cn)$. 
Our result strengthens Holmsen's bound by improving the dependence on $c$ from $c^{2^{r-1}}$ to linear in $c$ with bounded number of forbidden structures. 
Furthermore, our approach is naturally linked to the notion of VC-dimension.
\end{abstract}

\section{Introduction}
%==================motivations==================
Let $G$ be a graph with vertex set $V(G)$ and edge set $E(G)$. Denote by $\chi(G)$ the \emph{chromatic number} of $G$, and by $\omega(G)$ the \emph{clique number} of $G$, which is the order of the largest clique in $G$. Given a graph $H$, a graph $G$ is said to be (\emph{induced}) $H$-\emph{free} if it contains no (induced) subgraph isomorphic to $H$. 
The Tur\'an number ex$(n, H)$ denotes the maximum number of edges in an $H$-free graph on $n$ vertices. The classical case where $H = K_k$, the complete graph on $k$ vertices, was resolved by Tur\'an~\cite{turan1941}, who also characterized the unique extremal graph achieving this bound.  
For a general graph $H$, the Erd\H{o}s-Stone-Simonovits theorem (ESS-theorem)~\cite{erdos1966limit,erdos1946linear}  states that
$$
 \text{ex}(n, H) = \left(1 - \frac{1}{\chi(H) - 1}\right) \binom{n}{2} + o(n^2).
$$
By the Erd\H{o}s--Stone--Simonovits theorem, forbidding a fixed (non-induced) subgraph $H$ forces extremal $H$-free graphs to be asymptotically $(\chi(H)-1)$-partite (for $\chi(H)\ge3$), and hence their clique number is bounded by a constant depending only on $H$. 
While $\chi(H)=2$ (i.e. $H$ is bipartite) the theorem only yields $\mathrm{ex}(n,H)=o(n^2)$, which does not determine the clique number in general, but it does rule out cliques of linear size. In contrast, such a phenomenon completely disappears when $H$ is forbidden as an induced subgraph. 
In 2002, Gy\'arf\'as, Hubenko, and Solymosi~\cite{gyarfas2002} proved that if an $n$-vertex graph $G$ has at least $c \binom{n}{2}$ edges and contains no induced $K_{2,2}$, then $\omega(G)\geq \frac{c^2}{10} n$. 
This result was later generalized by Loh et al.~\cite{loh2016} to graphs that are induced $K_{2,t}$-free. 
More recently, Holmsen~\cite{holmsen2020} improved the result of Gy\'arf\'as et al., who showed that excluding an induced $K_{2,2}$ guarantees $\omega(G) \geq (1 - \sqrt{1 - c})^2 n$.

\begin{theorem}[\cite{holmsen2020}]\label{thm:k22induced}
Let $G$ be a graph on $n$ vertices with at least $c \binom{n}{2}$ edges. If $G$ contains no induced $K_{2,2}$, then
\[
\omega(G) \geq (1 - \sqrt{1 - c})^2 n.
\]
\end{theorem}
When $c$ is small, one finds that the bound of clique number in Theorem \ref{thm:k22induced} is $\Omega(c^2n)$, which is optimal as demonstrated by the blow-up of a projective plane graph. 
By forbidding more induced substructures, one can obtain better quantitative bound. A graph \(G\) is called \emph{chordal} if every cycle of length at least four in \(G\) has a chord, that is, \(G\) contains no induced cycle of length at least four.
Abbott and Katchalski~\cite{abbott1979}, and independently Gy\'arf\'as, Hubenko, and Solymosi~\cite{gyarfas2002}, established a corresponding tight bound for chordal graphs, improving the dependence on \(c\) from \(c^2\) to \(c\).

\begin{theorem} [\cite{abbott1979,gyarfas2002}]\label{thm:chor}
Suppose that $G$ is a chordal graph on $n$ vertices with at least $c \binom{n}{2}$ edges, then
\[
\omega(G) \geq (1 - \sqrt{1 - c}) n.
\]
Moreover, the chordal graph which is the sum of an independent set of order $t=\lceil\sqrt{1 - c}n\rceil$ and a complete graph $K_{n-t}$ shows that the inequality is asymptotically best possible.
\end{theorem}
Note that the bound of clique number in Theorem \ref{thm:chor} is $\Omega(cn)$ when $c$ is sufficiently small. Holmsen~\cite{holmsen2020} also extended the case $K_{2,2}$ to that of the multipartite version, see Theorem \ref{thm:r-blow2}. He in fact proved a more general hypergraph version, forbidding certain family of induced substructures. 
Using his result, one can prove that abstract colorful Helly implies abstract fractional Helly. The $t$-\emph{blow-up} of $G$, denoted by $G[t]$, is the graph obtained by replacing each vertex of $G$ with an independent set of size $t$, and each edge of $G$ with a complete bipartite graph between the corresponding sets. 

\begin{theorem}[\cite{holmsen2020}]\label{thm:r-blow2}
Let $G$ be a graph on $n$ vertices with at least $c \binom{n}{r}$ $r$-cliques. If $G$ contains no induced $K_r[2]$, the $2$-blow-up of a clique $K_r$, then there exists $M>0$ such that 
\[
\omega(G) \geq Mc^{2^{r-1}} n.
\]
\end{theorem}

Although this result is linear in $n$, its bound on $c$ deteriorates rapidly as \(r\) grows, making it less effective for large \(r\), whether this bound can be improved remains open. It is natural to ask whether Theorem~\ref{thm:r-blow2} can be strengthened in a manner analogous to the improvement from Theorem~\ref{thm:chor} to Theorem~\ref{thm:k22induced}. 
This forms the main result of this paper.

%extending the chordal graph setting to a multipartite generalization, which is the main focus of this paper.

\begin{definition}
Let $u_i$ and $u_i'$ be the two vertices in the $i$-th part of $K_r[2]$ and $U' = \{ u_1', \dots, u_r' \}$.  Define $\mathcal{K}_r^{[2]}$ as the family of subgraphs of $K_r[2]$ whose edge set is of the form $(E(K_r[2]) \setminus E(K_{U'})) \cup E'$, for some $E' \subseteq E(K_{U'})$, where $K_{U'}$ is the subgraph of $K_r[2]$ induced by $U'$.
\end{definition}

\begin{figure}[htbp]
    \centering
    \begin{tikzpicture}[scale=1.2,
      dot/.style={circle, fill=black, inner sep=1.5pt},
      every label/.style={font=\small}
    ]

    % 顶点位置
    \node[dot, label=left:$u_1$] (u1) at (0,0) {};
    \node[dot, label=left:$u_1'$] (u1p) at (0,-0.7) {};

    \node[dot, label=right:$u_2$] (u2) at (3,0) {};
    \node[dot, label=right:$u_2'$] (u2p) at (3,-0.7) {};

    \node[dot, label=left:$u_3$] (u3) at (1.15,1.6) {};
    \node[dot, label=right:$u_3'$] (u3p) at (1.85,1.6) {};

    % 黑色边（非 U' 内部边）
    \foreach \i/\j in {u1/u2, u1/u2p, u1p/u2, u1p/u2p,
                       u1/u3, u1/u3p, u1p/u3, u1p/u3p,
                       u2/u3, u2/u3p, u2p/u3, u2p/u3p} {
        \draw[black] (\i) -- (\j);
    }

    % 红色边：U' 内部完全图
    \draw[red, thick, line cap=round] (u1p) -- (u2p);
    \draw[red, thick, line cap=round] (u1p) -- (u3p);
    \draw[red, thick, line cap=round] (u2p) -- (u3p);

    \end{tikzpicture}
    \caption{The family $\mathcal{K}_3^{[2]}$ includes all subgraphs of $K_3[2]$ that contain all black edges and any subset of the red edges.}
    \label{fig:Kr-2-example}
\end{figure}
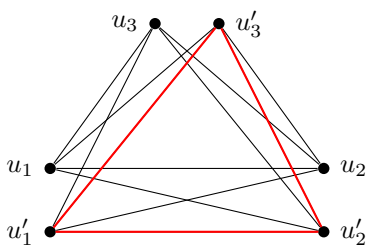

We give an example for $r=3$ to illustrate this notation in Figure \ref{fig:Kr-2-example}.

%example

%=====remark=========
%\paragraph{Remark.}
%For $r = 2$, the class of induced $\mathcal{K}_2^{[2]}$-free graphs---that is, graphs admitting no induced subgraph isomorphic to $P_4$ or $K_{2,2}$---is intimately related to the class of chordal graphs, although the two families do not coincide. Every induced $\mathcal{K}_2^{[2]}$-free graph is necessarily chordal, since the presence of any induced cycle of length at least four would yield an induced $C_4$ or $P_4$, both of which are prohibited substructures. Conversely, the inclusion is strict: not every chordal graph is induced $\mathcal{K}_2^{[2]}$-free. In particular, although chordal graphs are, by definition, free of induced cycles $C_\ell$ with $\ell \ge 4$, some chordal graphs---most notably the path $P_4$---still contain an induced $P_4$ and hence fail to satisfy the $\mathcal{K}_2^{[2]}$-free property. It follows that while every induced $\mathcal{K}_2^{[2]}$-free graph is chordal, the converse holds only for those chordal graphs devoid of induced $P_4$'s. In particular, a forest is chordal, but it is induced $\mathcal{K}_2^{[2]}$-free if and only if it contains no induced $P_4$. Consequently, among all chordal graphs, the unique obstruction to being induced $\mathcal{K}_2^{[2]}$-free is the path $P_4$ itself.
\paragraph{Remark.}
    For $r = 2$, the class of induced $\mathcal{K}_r^{[2]}$-free graphs---those containing no induced $P_4$ or $K_{2,2}$---resembles the chordal graphs. It forms a subclass of chordal graphs, since any induced cycle of length at least $4$ would contain an induced $P_4$ or $K_{2,2}$. However, there is a subtle difference. While  every chordal graph is induced $K_{2,2}$-free, some, such as $P_4$, still contain induced $P_4$'s and hence are not induced $\mathcal{K}_2^{[2]}$-free.

\vspace{0.3cm}
Motivated by the role of chordal graphs in the case $r=2$, we extend the approach to general $r$ by considering induced $\mathcal{K}_r^{[2]}$-free graphs, which serve as a natural multipartite analogue in this setting. This extension leads to an improved bound that is linear in $c$.

\begin{theorem}\label{main}
 Let $r \ge 2$ and $0<c<1$, and let $G$ be a graph on $n$ vertices with at least $c \binom{n}{r}$ copies of $K_r$. If $G$ is induced $\mathcal{K}_r^{[2]}$-free, then, for sufficiently large $n$, $\omega(G) \ge \frac{c}{18r}n$. 
\end{theorem}

Our result concerns graphs that forbid bounded number of  induced substructures (at most $2^{\binom{r}{2}}$). 
For \(r=2\), this family corresponds to the induced \(\mathcal{K}_2^{[2]}\)-free graphs, which exclude just $2$ induced configurations, in contrast to chordal graphs that forbid infinitely many. 
Interestingly, our theorem attains a bound of the same order as that for chordal graphs, showing that far fewer structural restrictions already suffice to ensure linear-size cliques.

\vspace{0.4cm}
\noindent
\textbf{Related work:} Our proof strategy relies on the notion of VC-dimension, whose connection to bi-induced subgraphs has been established in earlier work  \cite{nguyen2023}, with Example \ref{eg:bi-in} providing a concrete illustration of this correspondence. Our extension to semi-induced subgraphs (see Claim \ref{claim:VC-dimension}) constitutes a natural generalization of the bi-induced case and yields a broader understanding of the role of VC-dimension in graph structures.

%==================================
%-----------------------------before----------------

\section{Proof of Theorem~\ref{main}}
Before proving Theorem~\ref{main}, we first introduce the following definitions.
Let $\mathcal{F}$ be a set system on a ground set $X$. 
For any subset $S\subseteq X$, define $\mathcal{F}|_{S} := \{ F \cap S: F \in \mathcal{F} \}$. 
A set $A\subseteq X$ is $shattered$ by $F$ if $\mathcal{F}|_{A} = 2^{A}$, that is, for every $A'\subseteq A$ there exists $F\in \mathcal{F}$ such that $F \cap A = A'$.
The \emph{VC}-\emph{dimension} of $\mathcal{F}$ is the maximum size of a set $A\subseteq X$ that is shattered by $\mathcal{F}.$ 

%=============adding==========
For general graph $G$, the VC-dimension of the family of neighbourhood of vertices in $G$: 
\[
\mathcal{F}_{G} = \{N_G(v) : v \in V(G)\}
\] 
is often taken as the VC-dimension of the graph $G$ itself, as discussed in~\cite{nguyen2023}.  
A set $A \subseteq V(G)$ of size $t$ shattered by $\mathcal{F}_G$ can be interpreted as a particular structure within $G$ as follows. 

\begin{example}\label{eg:bi-in}
\rm Let $A$ be a $t$-subset of $V(G)$, and let $B$ be a set of $2^t$ vertices disjoint from $A$.  
Consider a bijection $\sigma : B \to 2^A$ from $B$ to the power set of $A$.  
Define the graph $F_{A,B;\sigma}$ to be the subgraph induced by $A \cup B$ such that: for edges between $A$ and $B$, each vertex $v \in B$ is adjacent exactly to the vertices in $\sigma(v)$,
and for edges inside $A$ and inside $B$, any subgraph can be  permitted (e.g., arbitrary edges).    
\end{example}

If $A$ is shattered by $\mathcal{F}_G$, then there exist a set $B$ and a bijection $\sigma$ such that $F_{A, B;\sigma}$ is an induced subgraph of $G$. When considering the VC-dimension of the set of maximal cliques in $G$, a shattered set $A$ of size $t$ corresponds to an induced subgraph $F$ of $G$ with $F \in \mathcal{K}_t^{[2]}$.

Next we introduce the following lemma, illustrating that bounded VC-dimension restricts a set family to polynomially many patterns on any finite set, which is essential to the proof of Theorem \ref{main}.

\begin{lemma}[Sauer-Shelah~\cite{sauer1972,shelah1972}]\label{SS}
    For any $\mathcal{F}$ with VC-dimension $k$ and for any $|S|=m$, we have $$\big|\mathcal{F}|_{S}\big|\leq \binom{m}{0}+\binom{m}{1}+\cdots+\binom{m}{k}.$$
\end{lemma}

\begin{proof}[Proof of Theorem ~\ref{main}]
Let $G$ be an induced $\mathcal{K}_r^{[2]}$-free graph on  $n$ vertices containing at least $c\binom{n}{r}$ copies of $K_r$. In the following, we will prove that $G$ contains a clique of size $\frac{c}{18r}n$. 

Let $\mathsf{MC}(G)$ be the set of maximal cliques in $G$. 
%Let $MIS(G)$ be the set of all maximal independent sets in $G$.
%Let $\mathcal{F}$ be the dual class of $\mathsf{MC}(G)$, that is, $$\mathcal{F} := \{ E_v: v \in V(G) \}, \  E_v := \{ K\in \mathsf{MC}(G): v\in K \}.$$  %\coloneqq

\begin{claim}\label{claim:VC-dimension}
    The VC-dimension of $\mathsf{MC}(G)$ is at most $r-1$.
\end{claim}
\begin{proof}
    Suppose, to the contrary, that there exists a set $S \subseteq V(G)$ of size $r$ such that $\mathsf{MC}(G)|_{S} = 2^{S}$, and write $S = \{ u_1, \dots, u_r \}$. For each $i \in [r]$, let $S_i = S \setminus \{ u_i \}$. Then there exists a maximal clique $K^i \in \mathsf{MC}(G)$ such that $K^i \cap S = S_i$; 
    that is, $u_j \in K^i$ for all $j \neq i$ and $u_i \notin K^i$. This implies that $S$ forms a clique in $G$. By the maximality of $K^i$, there exists a vertex $u_i' \in V(K^i)\setminus S$ such that $u_i u_i' \notin E(G)$. 
    Consequently, for any distinct $i, j \in [r]$, we have $u_i u_j, u_i' u_j \in E(G)$, and $u_i u_i' \notin E(G)$. 
    Hence, the subgraph induced by the set $\{ u_1, u_1', \dots, u_r, u_r' \}$ is isomorphic to a graph in $\mathcal{K}_r^{[2]}$, which is a contradiction.
\end{proof}

We fix $m = \left\lfloor \frac{9r}{c} \right\rfloor$ and let $c' = \frac{c}{18r}$. Then
$$\sum_{i=0}^{r-1} \binom{m}{i} \le 2\binom{m}{r-1} < \frac{1}{4}c \binom{m}{r}.$$
Here, we also assume that $n$ is sufficiently large so that $c'n \geq m$.

For the sake of contradiction, suppose that $\omega(G) < c'n$. Let $S_m \subseteq V(G)$ be a subset of size $m$. By Lemma~\ref{SS} and Claim \ref{claim:VC-dimension}, we have  
\begin{equation}\label{equ:MC}
\left| \mathsf{MC}(G)|_{S_m} \right| \le \sum_{i=0}^{r-1} \binom{m}{i} < \frac{1}{4}c \binom{m}{r}.
\end{equation}

Now we choose a pair $(S_r, S_m)$ uniformly at random with $S_r \subseteq S_m$, and examine the probability that there exists a maximal clique $K$ in $G$ satisfying $S_r = S_m \cap K$. Note that any set $S_r$ satisfying the above condition belongs to $\mathsf{MC}(G)|_{S_m}$.

%It follows that the probability that $S_m$ contains $S_r$ but $S_m \setminus S_r$ contains no vertex of $K$ is at least

Since $G$ contains at least $c\binom{n}{r}$ copies of $K_r$, the probability that a randomly chosen $S_r$ induces a clique in $G$ is at least $c$. Suppose that $S_r$ indeed induces a clique in $G$. By our assumption, this clique is contained in a maximal clique $K$ of size less than $c'n$. 
It follows that the probability that $S_m$
contains $S_r$ but $S_m
\setminus S_r$ contains no vertex of $K$ is at least
\begin{equation}\label{equ:pi}
\frac{\binom{n - c'n}{m - r}}{\binom{n - r}{m - r}} \ge \prod_{i=0}^{m - r - 1} \frac{n - c'n - i}{n - i} \ge \left( \frac{n - c'n - m}{n - m} \right)^m\geq(1 - 2c')^m \ge e^{-2c'm} \geq \frac{1}{4},
\end{equation}
where the third-to-last inequality holds as $c'n \ge m$ and recall that $c' \leq \frac{1}{2m}$.

Hence, the probability that the pair $(S_r, S_m)$ satisfies $S_r = S_m \cap K$ for some maximal clique $K$ in $G$ is at least $\frac{1}{4}c$.

For a random set $S_m$, we find that the expected number of sets $S_r$ such that $S_r = S_m \cap K$ for some maximal clique $K$ in $G$ is at least
$$\frac{1}{4}c \binom{m}{r}.$$ 
Thus, there exists some $S_m$ such that
$$\left| \mathsf{MC}(G)|_{S_m} \right| \ge \frac{1}{4}c \binom{m}{r},$$
which contradicts to (\ref{equ:MC}), and then we have $\omega(G) \geq c'n=\frac{c}{18r}n$. 
\end{proof}

\section{Concluding remarks}
In this paper, we improve the bound from $\Omega(c^{2^{r-1}})$ to $\Omega(c/r)$ by forbidding $2^{\binom{r}{2}}$ graph structures. It would be interesting to know when the transition of $c^{2^{r-1}}$ to linear in $c$ happens. More precisely, what is the minimum number of forbidden graph structures that forces a bound linear in $c$?

\vspace{0.5cm}
\textbf{Acknowledgement.}
The work in this article was carried out in the ECOPRO group at IBS. We are grateful to Professor Hong Liu for suggesting this research direction.

\bibliographystyle{abbrv}
\bibliography{reference}

\end{document}